\documentclass[11pt]{article}

\def \S{{\bf /\!\!\!\! S}}
\def \DD {{\bf /\!\!\!\!D}}
\def \C{{\bf C}}
\def \Cl{{\bf C l}}
\def\I{{\cal I}}
\def\e{{\varepsilon}}
\def\oe{{{\overline\varepsilon}}}
\def\oc{{\cal C}}
\def \Z{{\bf Z}}
\def \R{{\bf R}}

\def\H{{\cal H}}
\def\V{{\cal V}}
\def\P{{\bf P}}

\newcommand{\Kahler}{K\"{a}hler}

\newtheorem{defn}{Definition}[section]

\newtheorem{lemma}[defn]{Lemma}

\newtheorem{theorem}[defn]{Theorem}
\newtheorem{definition}[defn]{Definition}
\newtheorem{remark}[defn]{Remark}
\newtheorem{proposition}[defn]{Proposition}
\newtheorem{corollary}[defn]{Corollary}

\newcommand{\bg}{{\bf
 \bigtriangledown}}

\begin{document}
\title{Anti-holomorphic twistor and symplectic structure
\thanks{{\it Key words and phases.} twistor space, pure spinor, symplectic structure}}

\author{Dosang Joe\thanks{ This work was supported by the BK21 project.}
\\Ewha Women's university}
\date{May 26, 2000}

\maketitle
\begin{abstract}
It is well known that the twistors, section of twistor  space,
classify the almost complex structure on even dimensional
Riemannian manifold $X$. We will show that existence of a harmonic
and  anti-holomorphic  twistor  is equivalent to having a
symplectic structure on $X$
\end{abstract}

\section{Introduction}
Recently, the interest  of symplectic manifolds has been growing
in a perspective of Mathematical Physics related field,   for
example, Quantum cohomology theory, Seiberg-Witten theory etc. By
definition, a manifold having  a non-degenerate closed two form
$\omega$ is called  a symplectic manifold. This category of
manifolds was firstly understood as that  of \Kahler \, \,
manifolds, which has even odd betti number, for example, later on
some mathematician like B. Thurston and R. Gompf  constructed
examples of symplectic manifolds which cannot have \Kahler\,\,
structure. Moreover R. Gompf \cite{Go} find a systematic way of
constructing symplectic manifolds and show that every finitely
presented group can be realized as a fundamental group of a
symplectic 4-manifold. It reveals that the symplectic category is
much more bigger than \Kahler\,\, one and expected to be
characterized as cohomology condition of given manifold such as
$a\in H^2 (X, \R) $ and $0 \ne a\cup \cdots \cup a \in H^{2m} (X,
\R)$. This expectation has been broken in advent of Seiberg-Witten
theory for the 4-dimensional topology. It has been known that
every symplectic 4-manifold has non-zero Seiberg-Witten invariants
\cite {T1, T2}, which indicates that condition of having
symplectic structure on 4-manifolds is quite subtle. Taking closer
look at the Taubes's paper \cite{T1}, we can find that he was
making use of the characterization of symplectic form, which is
there are canonical $Spin^c$ structure associated almost complex
structure $J$ and naturally induced a nowhere vanishing positive
spinor $u$ which is harmonic ,i.e., $\DD u=0$.  In this paper, we
are going to show that this characterization is equivalent to the
existence of the symplectic form on a given manifold.  First of
all, notice that symplectic form $\omega$  on a given manifold
realized as an imaginary part of Hermitian metric for some almost
complex structure $J$ on $TX$. Hence the existence of almost
complex structure is a necessary condition for that of symplectic
structure. Given  a Riemannian even dimensional manifold, the
orthogonal almost complex structure is equivalent to  a section of
the twistor space which is a canonical fiber bundle of
$SO(2m)/U(m)$. We will discuss on this in the section 2. After
choosing a twistor $u$ , equivalently having an almost complex
structure $J$), there is the naturally associated $Spin_2m$
equivariant  Hermitian metric on $(TX, J)$ and a canonical
$Spin^C$ representation. The imaginary part of the Hermitian
metric $\omega$ is our candidate for the symplectic form. It is
easily derived that the condition for $d\omega =0$ is equivalent
to the section $u$ is anti-holomorphic and harmonic ($\DD u=0$),
where $u$ can be understood as a positive spinor of the canonical
positive spinor bundle.  To prove this theorem is main purpose of
this paper. It also gives a simple characterization of symplectic
structure on smooth 4-manifolds, which is the same as the Taubes'
analysis of symplectic form. Conclusively, the number of equations
for $\omega $ being a symplectic form is $m(m-1)/2 +
m(m-1)(m-2)/6$ which is bigger than that for integrability
condition which is $m(m-1)/2$. That is a little bit odd because
the space of symplectic form (deformation space; {\it it is open
in the space of two form $\Omega^2(X)$})  is rather larger than
the integrable complex structure which is finite dimension. On the
other hand, since the deformation space is a kind of big,  there
are a lots of such an anti-holomorphic and harmonic twistor for
some Riemannian metric on $TX$ once $X$ has a symplectic
structure. It gives rise a question whether  the condition we have
found is ``generic", which means in symplectic manifold, the
generic Riemannian metric can be induced by a symplectic form
$\omega$ and
 an almost-complex structure $J$ associated to it.
\section{ Pure Spinor and twistor}
Fix $\R^n$ be the standard inner product( $<,>$)  real vector
space and extend this metric $\C$-linearly to $\C^n =R^n
\otimes\C$. Let  $\Cl_n =Cl_n\otimes \C$ be the associated
complexified Clifford algebra. Let $\S_\C$ be the fundamental
$\Cl_n$-module which defines the irreducible complex spinor space.
For each spinor $\sigma \in \S_\C $, we can consider the
$\C$-linear map
 $$
j_\sigma : \C^n \to \S_\C \quad \mbox{ given by} j_\sigma
(v)\equiv v\cdot \sigma $$ Generically, this map is injective.
However, there are interesting spinors for which $ \dim(\ker
j_\sigma)
>0$.
\begin{defn}
A complex subspace $V \subset \C^n$ is said to be {\bf isotropic}
(with respect to the bilinear form $<\cdot , \cdot>$) if $<v,
w>=0$ for all $v, w \in V$.
\end{defn}

We define a hermitian inner product $(\cdot, \cdot)$ on $\C^n$ by
setting $(v,w)= <v,\overline{w}>$. Clearly, if $V \subset \C^n$ is
an isotropic subspace, then $V\perp\overline{V}$ in this hermitian
inner product. In particular, therefore, we have
$$
 2\dim_\C V \le n.
$$
\begin{definition}
A spinor $\sigma$ is {\bf pure} if $\ker j_\sigma$ is a maximal
isotropc subspace, i.e., if $\dim(\ker j_\sigma)= \mbox{[n/2]}$.
\end{definition}
Denote by $P\S$ the subset of pure spinors in $\S_C$, and denote
by $\I_n$ the set of maximal isotropic subspaces of $\C^n$ Both
$P\S$ and $I_n$ are naturally acted upon by the group $Pin_n$, and
the assignment $\sigma \mapsto \ker j_\sigma$ gives a
$Pin_n$-equivariant map
 $$
 K : P\S \longmapsto \I_n.
 $$
From this point on we shall assume that $n=2m$ is an even integer,
and furthermore that $\R^{2m}$ is oriented.
 \begin{definition}
 An {\bf orthogonal almost complex structure} on $\R^{2m}$ is an
 orthogonal transformation  $J : \R^{2m}\to R^{2m} $ which satisfies
 $J^2=-\mbox{Id}$. For any such $J$, an associated {\bf unitary
 basis} of $\R^{2m}$ is an ordered orthonormal basis of the form
$ \{e_1, Je_1, \cdots e_m, Je_m\}$. Any two unitary bases for a
 given  $J$ determine the same orientation. This is called the
 {\bf  canonical orientation} associated $J$.
\end{definition}
 Let $\oc_m$ denote the set of all orthogonal almost complex
 structures on $\R^{2m}$. It is easy to see that $\oc_m$ is a
 homogeneous space for the group $O_{2m}$.  It falls into two
 connected components $\oc_m^+$ and $\oc_m^-$ where $\oc_m^+ \cong
 SO_{2m} /U_m$ consists of those almost complex structures whose
 canonical is {\bf positive} ( i.e. agrees with given one on
 $\R^{2m}$). Associated to any $J\in \oc_m$ there is a decomposition
$$ \C^{2m} =V(J)\oplus \overline{V(J)}, \mbox {where} $$
$$V(J)\equiv \{ v\in \C^{2m} : Jv=-i v\} = \{v_0 +iJv_0 : v_0 \in
\R^{2m}\}
 $$

There is an $O_{2m}$-equivalent bijection $$ \oc_m
\stackrel{V}{\longrightarrow}\I_{2m}$$ which associates to $J$ the
isotropic subspace $V(J)$ Let $\I^+_{2m}$ denote the component
corresponding  to $\oc_m^+$. Using the complex volume element
$\omega_\C = i^m e_1\cdots e_{2m}$, we have  a decomposition
$\S_\C =\S_\C^+ \oplus \S_\C^-$ into $+1$ and $-1$ eigenspace
respectively. Easy calculation gives  a decomposition $P\S =P\S^+
\coprod P\S^-$ of the pure spinor space into positive and negative
types. Let $\P (P\S^+)$ denote the projectivization of the pure
spinor space, i.e., $\P (P\S^+)=P\S/\sim$ where we say that
$\sigma\sim\sigma^\prime$ if $\sigma= t\sigma^\prime$ for some
$t\in \C$. Each of the space $\P (P\S^\pm), \oc^\pm_m$ and
$\I^\pm_{2m}$ are acted upon by $Spin_{2m}$, in fact by $SO_{2m}$.
\begin{proposition}
The maps $\sigma \mapsto K(\sigma)$ and $J\mapsto V(J)$ induce
$SO_{2m}$-equivariant diffeomorphisms $$ \P
(P\S^+)\stackrel{K}{\mapsto}\I^+_{2m} \stackrel
{V}{\mapsto}\oc^+_m \quad\quad \mbox{ and}\quad \quad\P
(P\S^-)\stackrel{K}{\mapsto}\I^-_{2m} \stackrel
{V}{\mapsto}\oc^-_m $$
\end{proposition}
We refer to the original book \cite{LM} for details.

For the sake of further discussion,  we will fix $V\in I_{2m}^+$
and let $J\in \oc_m^+$ be the associated complex structure. Choose
a unitary basis $\{ e_1, Je_1, \cdots, e_m, Je_m\}$ of $\R^{2m}$
and set $$\e_j ={1\over \sqrt{2}} (e_j -iJe_j) \quad\quad
\overline{\e}_j={1\over \sqrt{2}} (e_j +iJe_j).$$ Define
\begin{equation}
\omega_j=-\e_j\overline{\e}_j \quad\quad
\overline{\omega}_j=-\overline{\e}_j\e_j \end{equation} Let $W$ be
a linear subspace invariant under multiplication by $e_j$ and
$Je_j$. Then there is a hermitian orthogonal direct sum
decomposition $$W=W_j \oplus W_j^\prime$$ where $$
W_j=\overline{\omega}_j\cdot W=\ker(\mu_{\overline {\e}_j }|_W)
\quad \mbox{and} \quad W_j={\omega}_j\cdot W=\ker(\mu_ {\e_j
}|_W)$$ and where $\mu_{\e_j} : W\to W$ is defined by $\mu_{\e_j}
(w)= \e_j\cdot w.$ By direct inductive calculation, we can
construct $$\S_m =\ker(\mu_{\overline{\e}_1}) \cap \cdots \cap
\ker(\mu_{\overline{\e}_m}) \quad \dim_\C \S_m= 1$$ The complex
volume form $\omega_\C =i^m e_1 Je_1 \cdots e_m Je_m$ has the
value $+1$ on $\S_m$ because $\overline{\e}_j \sigma =0
\Rightarrow -ie_jJe_j \sigma =\sigma.$ Therefore, $\S_m \subset
\S_\C^+$.

We clearly have that $V(J)=\ker j_\sigma $ for $\sigma\in \S_m$.
Hence $\S_m$ is independent of the choice of unitary basis and the
map $V \mapsto [\S_m]$ gives the desired map $K^{-1}$ for the
above proposition.
\begin{definition} The bundle $\tau(X) \cong \P
(P\S^+)$ is called the {\bf twistor space} of $X$.
\end{definition}
Note that $\P (P\S^+)$ is an $SO_{2m}$-bundle and is globally
defined whether or not $X$ is a spin manifold.

The total space of $\tau(X)$ carries a canonical almost complex
structure defined by using the canonical decomposition of tangent
space of $\tau(X)$, which is induced by the Riemannian connection
of $X$. $$ T(\tau(X))= \V \oplus \H$$ where $\H$ is a field of
horizontal planes and $\V$ is the field of tangent planes to the
fibers.  As noted, $\V$ has an almost complex structure integrable
on the fibers since the fiber is naturally homogeneous complex
manifold ($\cong SO_{2m}/U(m)$). The bundle $\H$ has a
``tautological" almost complex structure defined, via the
identification $\pi_\ast : \H_J \to TX$, to be the structure $J$
itself.

The question of integrability of $J$ already accomplished by M.
Michelsohn.
\begin{theorem}
{\rm \cite {LM, M}} Let $X$ be an oriented(even-dimensional)
riemannian manifold with an almost complex structure determined by
a projective spinor field $u\in \Gamma(\tau(X))$. Then this almost
complex structure is integrable if and only if $u$ is holomorphic.
\end{theorem}
This will be proved in Remark 3.3.  As mentioned above, $\tau(X)$
carries a canonical almost complex structure. Now a $C^1$-map
between almost complex manifolds $f : (X, J_X) \to (Y, J_Y)$ will
be called {\bf holomorphic (resp. anti-holomorphic)} if its
differential $f_\ast$ is everywhere $J$-linear(resp.
anti-$J$-linear) i.e., if $f_\ast\circ J_X =\pm J_Y \circ f_\ast$
respectively.
\begin{remark}
{\rm More succinctly one could say that {\it cross-section of
$\tau (X)$ induce almost complex structure, and holomorphic
cross-section induce the integrable ones} However, the condition
that a cross-section $u$ be holomorphic is {\it not} linear since
the complex structure on $X$ depends itself on $u$. }
\end{remark}
We will prove that the complimentary condition for the
holomorphicity, which is anti-holomorphic and harmonic is
equivalent to that $u$ induce a symplectic structure on $X$

\begin{definition}
$\omega\in \Omega^2(X)$ is a symplectic form if it is
nondegenerate  closed form. Moreover, $(X, \omega)$ is called a
symplectic structure on $X$.
\end{definition}
Given a twistor $u\in \P (P\S)$, there is naturally associated
nodegenerate  differential 2-form. It is induced by the hermitian
metric with respect to the almost complex structure $J$ and
Riemannian metric $g$ on $TX$ i.e., $$\omega (v,w) \equiv g(Jv,
w)$$ where $J$ is the almost complex structure corresponding to
$s\in \P (P\S)$. Moreover it can be written as in terms of unitary
basis, in other words,  $\omega =\sum_{i=0}^m e^\ast_i\wedge
(Je_i)^\ast$ where $e^\ast \in T^\ast X$ such that $e^\ast (v)=
g(e, v) \in \R$. Recall that $\omega_j =-\e_j \oe_j $ for complex
unitary basis $\{\e_1,\cdots ,\e_m , \overline{\e}_1 \cdots
\overline{\e}_m\}$ of $(TX\otimes \C)$. Since $\omega_j = -\e_j
\oe_j = 1-i e_j\cdot Je_j$, $i \omega = m -\sum_j \omega_j $.
\begin{eqnarray*}
{\omega}_1 \cdots {\omega}_m &=& \prod (1-ie_j\cdot Je_j )\\ &=&
1-i\sum e_j\cdot Je_j -\sum_{j\ne k} (e_j\cdot Je_j)\cdot
(e_k\cdot Je_k) +\cdots\\ &=& 1-i\omega +(1/2) (-1)^2
i\omega\wedge i\omega +\cdots + (1/m!)(-1)^m i\omega\wedge \cdots
i\omega\\ &=& 1-i\omega+(1/2!)(-i)^2\omega^2+\cdots +
(1/m!)(-i)^m\omega^m
\end{eqnarray*}
where $\omega^k =\overbrace{\omega\wedge \cdots \wedge
\omega}^{k\,\, times} \in \Omega^{2k}(X)$.

\begin{remark}
The above equality comes from the identification between $TX$ and
$TX^\ast$ via Riemannian metric. Note that $(1/m!) i^m \omega^m
=i^m e_1\cdot Je_1\cdots e_m\cdot Je_m= \omega_\C$
\end{remark}
Note that $\ast_\C \omega^k =k!/(m-k)! \omega^{m-k}$ i.e.,
$$d\omega =0 \Leftrightarrow  d\omega=d^\ast \omega=0
\Leftrightarrow \bigtriangleup_g ( {\omega}_1 + \cdots +{\omega}_m
)=0$$ where $\bigtriangleup_g$ is the Laplacian operator with
respect to metric $g$.  Hence we have that $\omega$ defines a
symplectic form if and only if $\tilde\omega = {\omega}_1 +\cdots
+{\omega}_m$ is harmonic. Our goal is to prove the following
theorem.
\begin{theorem}
 Let $X$ be an oriented(even-dimensional)
riemannian manifold with an almost complex structure determined by
a projective spinor field $u\in \Gamma(\tau(X))$. Then this almost
complex structure carries  symplectic structure  if and only if
$u$ is harmonic and anti-holomorphic.
\end{theorem}
 The
product element , $q=\overline{\omega}_1 \cdots
\overline{\omega}_m$ (conjugate of the above product ), of the
complexified Clifford algebra $\Cl_{2m}(X)$ can be  be
characterized at least locally by an element of $q \in End(\S^+)$
such that
\[ q ( \sigma )=\left\{ \begin{array}{clll} 0  & & \mbox{ if  }&
\sigma \in s^\perp \subset
 \S_\C \\ k\sigma \ & k\in \C^\ast  & \mbox {  and  if   }& [\sigma]=s \end{array}\right. \]
Note that we have not defined a  complex spin representation
$\S_\C$ globally over $X$. Without any specification of the
complex spinor bundle, the $\overline{\omega}_1 \cdots
\overline{\omega}_m $ is well-defined as an element of $\Cl_{2m}
(X)$. Using the almost complex structure associated with the
twistor $u$, we can define canonical $spin^c$ structure and
canonical complex spin representation. Given the canonical
$Spin^c$ representation, the product element $q = 2^m u\otimes
u^\ast$, which is an element of $q\in \mbox{End}_\C (\S^+)$ in a
way  of that $q(\alpha) = <\alpha, u > u $. In the next section,
we will prove that $(d\omega) u=0$ if and only if $\DD u=0$ by
using the action of $q$.
\section{$Spin^c$ representation and proof of Theorem 2.10}
Since $Spin^c_n \equiv Spin_n \times_{\Z_2} U(1)$, we have a short
exact sequence
 $$
 0\longrightarrow \Z_2 \longrightarrow Spin^c_n
 \stackrel{\xi}{\longrightarrow}
 SO_n\times U(1)\longrightarrow 1.
 $$ A principal $SO_n$-bundle $P$ carries a
$Spin^C$ structure if any only if the $w_2(P)$ is the mod 2
reduction of an integral class. Given  a twistor $u \in \P
(P\S^+)$, there is the  canonical  orthogonal almost complex
structure structure $J$ on $TX$  associated with $u$. This $J$
defines a canonical $Spin^c$ structure $\det_\C TX=K_X^{-1}$ since
the first Chern class of $K_X^{-1}$ is an integral lift of the
second Stiefel Whitney class, i.e., $c_1(K_X^{-1})\equiv w_2(X)
\mbox{mod 2}$. Let $\S_\C$ be the associated spinor bundle.  Using
the complex volume form $\sqrt{-1}^m e_1\cdot Je_1\cdots e_m
\cdots Je_m$, we have the decomposition of $\S_\C$ by the
$\pm$-eigenspace  of the complex volume element, where $\S^{\pm}=(
1\pm \omega_\C) \S_\C
$
Set $$\e_j ={1\over \sqrt{2}} (e_j -iJe_j) \quad\quad
\overline{\e}_j={1\over \sqrt{2}} (e_j +iJe_j).$$ be an unitary
basis for $TX$ as above. Define $$ \S_\C \cong \oplus
\S_{i_1,\cdots, i_m}\cong \oplus \ker (\mu_{\e_{i_1}}) \cap \cdots
\ker (\mu_{\e_{i_m}})
 $$
where  $\mu_{\e_{i_k}}= \left\{\begin{array}{cl} \mu_{\e_{k}} &
i_k= k \\ \mu_{\overline{\e_k}} & i_k=\overline{k}
\end{array}\right. $
Let $\sigma= \{i_1, \cdots i_m\}$ be the complex index used as
above, define $|\sigma|$ be the number of elements of the  subset
$\{ i_k =k\}$. Then we have $$ \S_\C^+ \cong
\displaystyle{\bigoplus_{|\sigma|=2i}} \S_{\i_1 \cdots i_m}\quad
\S_\C^- \cong \bigoplus_{|\sigma|=2i-1 } \S_{\i_1 \cdots i_m}$$
Especially, the twistor $u$ is contained in $ \S_{\overline{1},
\cdots \overline{m}}$ which is characterized as $\overline{\e}_j
\cdot u =0$ for all $j$. We can express the Dirac operator in
terms of the unitary basis, which follows that \begin{eqnarray*}
\DD& =& e_j \cdot \bg_{e_j} +Je_j \cdot \bg_{Je_j}\\
   &=& {1\over 2} (\e_j +\overline{\e}_j) \cdot \bg_{\e_j
   +\overline{\e}_j} - {1\over 2}(\e_j -\overline{\e}_j) \cdot \bg_{\e_j
   -\overline{\e}_j}\\
   &=& \overline{\e}_j\cdot \bg_{\e_j} +\e_j\cdot
   \bg_{\overline{\e}_j}
   \end{eqnarray*}

   \begin{remark}
  Note that  the covariant derivative $\bg$ is $Spin^c$ connection which is induced
  from both the Levi-Civita connection and the $U(1)$ connection
  on $K_X^{-1}$.  It should be well-noticed that our theorem is nothing to do with a $U(1)$
connection. Even though the condition we have imposed is related
to  simply ``local'' question, the $spin^c$ structure enable us to
work with globally.  Furthermore, the following argument we will
present below works finely without any $spin^c$ structure.
   \end{remark}

To define a Dirac operator on the spinors, we should specify a
$U(1)$ connection on $K_X^{-1}$. There is a canonical $U(1)$
connection unique up to gauge transformation $A_0$ such that $<\bg
u,  u>=0$. We will abuse  the notation $\DD$ for the Dirac
operator,  $\DD_{A_0}$,  which is induced by the Levi-Civita
connection and the canonical  $U(1)$ connection $A_0$.
  Our index notation convention  indicates that
$\bg_{\tilde e_j} \tilde e_k =\sum_l \tilde \omega_k^l (\tilde
e_j) e_l$ and $\Gamma_{j, k}^l =\tilde \omega_k^l (\tilde e_j)$ is
the Chistoffel symbol. Let $e_j = \tilde e_{2j-1}, Je_j =\tilde
e_{2j}$ then $\e_j = {1\over \sqrt{2}} (\tilde e_{2j-1} -i \tilde
e_{2j})$ where $i= \sqrt{-1}$. Then we have \begin{eqnarray*}
\bg_{\e_j} \overline{\e}_k =& a_{j,k}^{\,\,\, l}\e_l
+c_{j,k}^{\,\,\, l} \overline{\e}_l\,\, , \quad  \quad
\bg_{\overline{\e}_j} \e_k =& \overline{a}_{j,k}^{\,\,\,
l}\overline{\e}_l +\overline{c}_{j,k}^{\,\,\, l}
{\e}_l\\\bg_{\overline{\e}_j} \overline{\e}_k =& b_{j,k}^{\,\,\,
l}\e_l +d_{j,k}^{\,\,\, l} \overline{\e}_l\,\, , \quad  \quad
\bg_{{\e}_j} \e_k =& \overline{b}_{j,k}^{\,\,\, l}\overline{\e}_l
+\overline{d}_{j,k}^{\,\,\, l} {\e}_l
\end{eqnarray*} Since the Levi-Civita connection is naturally
compatible with the Hermitian metric on $TX\otimes \C$, we have
$$a_{j,k}^{\,\,\,l}= <\bg_{\e_j} \overline{\e}_k, \e_l>
=-<\overline{\e}_k,\bg_{\overline{\e}_j} \e_l> =- a_{j,l}^{\,\,\,
k}$$ By the same manner, we have $$ b_{j,k}^{\,\,\, l}
=-b_{j,l}^{\,\,\, k}\quad \mbox{ and} \quad c_{j,k}^{\,\,\, l} = -
\overline{d}_{j,l}^{\,\,\, k}$$ \begin{lemma} Let $u$ be a section
of twistor space and $J$ be the associated orthogonal almost
complex structure. Then $u$ is anti-holomorphic if and only if
$a_{j,k}^{\,\,\, l}=0$ for all $j, k, l$ and $u$ is holomophic
section if and only $b_{j,k}^{\,\,\, l}=0$ for all $j,k,l$.
\end{lemma}
 First of all, we have to
find the covariant derivative of $u$ which is $$ \bg u ={1\over 2}
\sum_{k < l} \tilde\omega^l_k\otimes \tilde e_l \tilde e_k \cdot u
$$ where $\tilde \omega$ is the $so(2m)$ connection 1-form(
Levi-Civita connection with respect to $g$) associated with
orthonomal basis $\{ \tilde e_1, \cdots \tilde e_{2m} \}$. Let  $$
\bg_{\e_j}u \equiv {1\over 2} \sum _{k < l} \tilde a_{j, k}^{\,
\,\, l} \otimes \e_l \cdot \e_k\cdot u\quad \mbox {mod } < u >$$
The coefficient $\tilde a_{j, k}^{\, \,\, l}$  can be derived as
follows, $$ \overline{\e}_t \cdot u = 0 \quad\quad\mbox {for all }
t$$ By taking covariant derivative $\bg_{\e_j}$, we have
$$(\bg_{\e_j} \overline{\e}_t) \cdot u + \overline \e_t \cdot
\bg_{\e_j} u =0$$ Hence
 \begin
{eqnarray*} (\bg_{\e_j} \overline {\e}_t ) \cdot u &=& -\overline
{\e}_t \cdot \sum_{k<l} {1\over 2}\, \tilde a_{j, k}^{\,\,\, l}
\e_l \e_k \cdot u\\ &=& -\sum {1\over 2}\, \tilde a_{j,k}^{\,\,\,
l} \, \overline{\e}_t \e_l \e_k \cdot u
\\ &=&\left \{ \begin{array}{cl}\,\tilde a_{j,k}^{\,\,\, l}\,
\e_l \cdot u &\mbox{ for  } k= t < l\\ -\,\tilde a_{j,k}^{\,\,\,t
} \e_k\, \cdot u & \mbox{ for } l=t> k
\end{array}\right.
\end{eqnarray*}
Since $\overline{\omega}_j=-\overline{\e}_j\cdot \e_j\cdot u =2
u$. Therefore $$ <\bg_{\e_j} \overline{\e}_t, \e_s> =  a_{j,
t}^{\,\,\, s }=\tilde  a_{j, t}^{\,\,\, s }$$ We get $\tilde
a_{j,k}^{\,\,\, l} = <\bg_{\e_j} \e_k, \e_l>.$ By the analogous
method, we can get $$\bg_{\overline{\e}_j}u \equiv {1\over 2} \sum
_{k < l} b_{j, k}^{\, \,\, l} \otimes \e_l \cdot \e_k\cdot u\quad
\mbox {mod } < u >.$$ With this understood, it can be rephrased
that $u$ is anti-holomorphic $\Leftrightarrow$ $\overline{\e}_t
\cdot \bg_{\e_j} u =0 $ $\Leftrightarrow$ $a_{j,k}^{\,\,\, l}
=0\Leftrightarrow <\bg_{\e_j} \overline{\e}_k , \e_l>=0$ for all
$j, k, l$.  Also  $u$ is holomorphic $\Leftrightarrow$
$\overline{\e_t} \cdot \bg_{\overline{\e}_j} u =0 $
$\Leftrightarrow$ $a_{j,k}^{\,\,\, l} =0\Leftrightarrow
<\bg_{\e_j} \overline{\e}_k , \e_l>=0$ for all $j, k, l$.
\begin{remark}  From the torsion free condition of Levi-Civita
connection, we have $$b_{j,k}^{\,\,\, l} -b_{k,j}^{\,\,\, l} =
<\bg_{\overline{\e}_j} \overline{\e}_k- \bg_{\overline{\e}_k}
\overline{\e}_j , \e_l> = < [ \overline{\e}_j, \overline{\e}_k ],
\e_l >.$$ Since the anti-commutativity  between upper index and
right lower index  $ b_{j,k}^{\,\,\, l}=-b_{j,l}^{\,\,\, k}$, we
can get an equivalent condition which says  that $b_{j,k}^{\,\,\,
l} -b_{k,j}^{\,\,\, l}=0$ if and only if $b_{j,k}^{\,\,\, l}=0.$
Hence it is easy to prove the Theorem 2.6 from the above equation.
\end{remark}
We want to find an equivalent condition for the harmonic two form
$\omega$ i.e. $\DD \omega =(d+d^\ast)\omega =0$, where $d^\ast$ is
the formal adjoint of $d$ with respect to $g$. The following lemma
is about it.
\begin{proposition} Let $\omega =-m+\sum_{k} {\omega}_k $
be the purely  imaginary part of the Hermitian metric. Then  $\DD
\omega =0$ if and only if $a_{j,k}^{\,\,\, l}=0$ and
$b_{j,k}^{\,\,\, l} +b_{l,j}^{\,\,\, k}+b_{k,l}^{\,\,\, j}=0$ for
all $j, k, l$.
\end{proposition} {\bf Proof:} Since $\omega$ is purely imaginary
two form, we have \begin{eqnarray*}\DD\omega &=&\sum_{j} \e_j
\bg_{{\oe}_j} \omega + \oe_j \bg_{\e_j}\omega
\\
&=& \sum_{j} \e_j \bg_{{\oe}_j} \omega -\overline{\sum_{j} \e_j
\bg_{{\oe}_j} \omega }\\ &=& 2i \mbox { Im }  \DD^{1\over2} \omega
\end{eqnarray*} It suffices to consider the half part of the
Dirac operator, it reads
 \begin{eqnarray*}-\DD^{1\over2} \omega &= &\sum_{j} -\e_j\cdot
 \bg_{{\overline\e}_j}\omega \\
&=&\sum_{j,k} ( \e_j\cdot (\bg_{\overline{\e}_j}\e_k)\cdot \oe_k
+ \e_j\cdot \e_k \cdot \bg_{\overline{\e}_j}\oe_k) \\ &=&
\sum_{j,k,l} (\overline{a}_{j,k}^{\,\,\, l} \e_j \oe_l \oe_k +
\overline{c}_{j,k}^{\,\,\, l} \e_j \e_l \oe_k +{b}_{j,k}^{\,\,\,
l} \e_j \e_k \e_l +
 {d}_{j,k}^{\,\,\, l} \e_j \e_k \oe_k )\\&=&
\sum_{j,k,l} (\overline{a}_{j,k}^{\,\,\, l} \e_j \oe_l \oe_k
+{b}_{j,k}^{\,\,\, l} \e_j \e_k \e_l )
+\sum_{j,k,l}(\overline{c}_{j,l}^{\,\,\, k}+{d}_{j,k}^{\,\,\, l})
\e_j \e_k \oe_l\\ &=&\sum_{i,j,l} (\overline{a}_{j,k}^{\,\,\, l}
\e_j \oe_l \oe_k +{b}_{j,k}^{\,\,\, l} \e_j \e_k \e_l )
\end{eqnarray*}
Hence $\DD \omega =0$ if and only if $\DD^{1\over2}\omega =0
\Leftrightarrow {a}_{j,k}^{\,\,\, l} =0 $ for all $i, j, k$  and $
\displaystyle{\sum_{\sigma}{b}_{\sigma(j)
,\sigma(k)}^{\,\,\,\,\,\,\, \sigma(l)}}$ where $\sigma$ is the
permutation of $i, j, k$. The relation  $b_{j,k}^{\,\,\,l} =-
b_{j,l}^{\,\,\,k}$ completes the proposition.

\noindent {\bf Proof of Theorem 2.10} Since $<\bg u, u>=0$, we
have
\begin{eqnarray*}\DD u&=& \sum_{j}(\oe_j\bg_{\e_j} u + \e_j
\bg_{\oe_j} u)\\ &=& \sum_{j,k,l}({1\over 4} a_{j,k}^{\,\,\, l}
\oe_j \e_k \e_l \cdot u +{1\over 4} b_{j,k}^{\,\,\, l}
\e_j\e_k\e_l \cdot u)
\\&=&\sum_{j,k} a_{j,j}^{\,\,\, k} \e_k \cdot u +\sum_{j<k<l}
{1\over 2}(b_{j, k}^{\,\,\, l}+ b_{l,j}^{\,\,\, k}+b_{k,l}^{\,\,\,
j}) \e_j\e_k\e_l\cdot  u
\end{eqnarray*}
Hence $u$ is anti-holomorphic pure spinor ($\oe_t \bg_{\e_j} =0$
for all $t,j$) and harmonic ($\DD u=0$ ) gives an equivalent
condition for $\omega$ being  a symplectic form. Note that given
symplectic manifold $(X,\omega)$ has such a anti-holomorphic and
harmonic twistor $u$ by choosing any almost complex structure
which calibrate $\omega$.

\begin{corollary}
$(d\omega) \cdot u=0 $  if and only if $u$ is harmonic , i.e.,
$\DD u=0$.
\end{corollary}
{\bf Proof:} Let $q=\prod_{j=1}^m \overline{\omega_j}=\prod_{j}
(1+ i e_j\cdot Je_j)$.  Using the action $q$ on $u$, $q\cdot u=
2^m u $, and taking Dirac operator on the both side, we can have
\begin{eqnarray*} \DD q\cdot u &=& (\DD q) \cdot u + \sum \tilde
e_j q\cdot \bg_{\tilde e_j} u\\
   &=& {\DD q} \cdot u \quad \Leftarrow (<\bg u, u> =0 \Rightarrow   q\cdot \bg_{\tilde e_j} u=0)\\
   &=& \DD u
   \end{eqnarray*}
   Thus  $\DD u= 0$ if and only if $(\DD q )\cdot u=0$. Moreover
since $(3i)^{m}(-1)^{{1\over 2} p(p+1)} \varphi \omega_{\C} = \ast
  \varphi$ for $\varphi \in \Omega^p (X)$ and $1/k! \ast\omega^k=1/(m-k)!
  \omega^{m-k}$, we have

\begin{eqnarray*} \DD q &=& (d + d^\ast) q = i d\omega +i^2/2! d\omega^2 +\cdots+
   i^{m-1}/(m-1)! d\omega^{m-1}\\
   && - i /(m-1)! \ast d \omega^{m-1}-i^2/(m-2)! \ast d
   \omega^{m-2} -\cdots  -i^{m-1} \ast d\omega\\
   &=&  i d\omega +i^2/2! d\omega^2 +\cdots+
   i^{m-1}/(m-1)! d\omega^{m-1}\\
   && i^{4m+1}  d\omega\cdot\omega_{\C}+i^{4m+2}/2! d\omega^2\cdot\omega_{\C} \cdots +i^{m-1}/(m-1)! d
   \omega^{m-1}\cdot\omega_{\C}
 \end{eqnarray*}
Since $\omega_{\C} u =u , \omega \cdot u=-(mi )u$, we have
$d\omega^k \cdot u = k (d\omega)\wedge \omega^{k-1} u =
k(-mi)^{k-1} d\omega \cdot u$. Thus $$(\DD q )\cdot u = 2i(1 +m+
m^2/2! +\cdots+m^{m-2}/(m-2)!) d\omega \cdot u .$$ This completes
the proof.
\begin{remark}
In dimension $2m\le 6$ every non-zero positive( or negative )
spinor is pure, i.e., $P \S^{\pm} = \S_{\C}^{\pm} -0$.  This is
simply because the group $Spin_{2m}$ acts transitively on the unit
sphere in $\S_{\C}^{\pm}$ in these dimensions.
\end{remark}

 In dimension 4, since $\varphi \in \Omega^3 (X, \R)$ acts on $u$
 injectively,
 we get $(d\omega)\cdot u=0$ if and only if $d\omega=0$
  Hence the harmonic spinor $u$, equivalently anti-holomorphic
  twistor, gives a sufficient condition to induce a symplectic structure.
  The next corollary follows from it.
\begin{corollary}
In dimension 4, Let $u$ be a nowhere vanishing section of positive
complex spinor bundle. Then $\DD u=0$ and $< \bg u, u>=0$ then $X$
is symplectic 4-manifold.
\end{corollary}
Finally, suppose  $\bg u=0$,  then  $u$ is then both holomorphic
and anti-holomorphic twistor. We have following corollary, which
is proposition 9.8 in \cite {LM}.
\begin{corollary}
Suppose $u$ is parallel, then $(X, g, J)$ becomes  a \Kahler\,\,
manifold.
\end{corollary}

\end{document}